\documentclass[12pt]{article}
\usepackage[latin1]{inputenc}
\usepackage{amssymb}
\usepackage{amsmath}
\usepackage{amsfonts}
\usepackage{amsthm}
\usepackage{enumerate}

\title{Singularity formation and blowup of complex-valued solutions of the modified KdV equation}
\author{J.L. Bona\thanks{Department of Mathematics, Statistics and Computer Science - University of Illinois at Chicago - Chicago, Il 60607 (U.S.A.)} 
\and S. Vento\thanks{Universit\'e Paris 13 - CNRS UMR 7539 Laboratoire Analyse, G\'eom\'etrie et Applications - 99 avenue J.B. Cl\'ement - 93430 Villetaneuse (FRANCE)} 
\and F.B. Weissler\footnotemark[2]}
\date{}

\numberwithin{equation}{section}

\newtheorem{theorem}{Theorem}[section]
\newtheorem{lemma}{Lemma}[section]
\newtheorem{proposition}{Proposition}[section]
\newtheorem{corollary}{Corollary}[section]
\newtheorem{definition}{Definition}[section]
\newtheorem{remark}{Remark}[section]
\newcommand\re[1]{(\ref{#1})}
\def\R{\mathbb{R}}
\def\Z{\mathbb{Z}}
\def\C{\mathbb{C}}

\def\N{\mathbb{N}}

\def\eps{\varepsilon}

\def\sech{\mathop{\rm sech}\nolimits}
\def\Im{\mathop{\rm Im}\nolimits}
\def\Re{\mathop{\rm Re}\nolimits}

\begin{document}
\maketitle
\noindent {\bf Abstract.}\,  The dynamics of the poles of the two--soliton solutions of 
the modified Korteweg--de Vries equation 
$$
u_t + 6u^2u_x + u_{xxx} = 0
$$
are determined.  A consequence of this study is the existence of classes of smooth, 
complex--valued solutions of this equation, defined for $-\infty <  x   < \infty$,  
exponentially decreasing to zero as $|x| \to \infty$, that blow up in 
finite time.  
\\

\section{Introduction}  Studied here is the modified Korteweg--de Vries equation 
\begin{equation}\label{mkdv1}
u_t + 6u^2u_x + u_{xxx} = 0,
\end{equation}
which has been derived as a rudimentary model for wave propagation in a number 
of different physical contexts.  The present paper is a sequel to the recent work
 \cite{BW2} wherein the 
dynamics of the complex singularities of the two--soliton solution
of the Korteweg de Vries equation,
\begin{equation}\label{kdv}
u_t + 6uu_x + u_{xxx} = 0,
\end{equation}
were examined in detail.  

The study of the pole dynamics of solutions of the Korteweg--de Vries 
equation and its near relatives began with some remarks of Kruskal 
\cite{Kr} in the early 1970's.   More comprehensive work was  carried 
out later, see \cite{Thi}, \cite{BowStu} and \cite{BryStu}.   One goal in our preceding paper \cite{BW2} was to
understand in more detail the propagation of  solitons in a neighborhood
of the interaction time.  Another motivation was an idea to be explained 
presently concerning singularity 
formation in nonlinear, dispersive
wave equations.  More particularly, we are interested in both  the generalized
Korteweg--de Vries equation
\begin{equation}\label{gkdv}
u_t + u_{xxx} + u^pu_x = 0
\end{equation}
and coupled systems of Korteweg--de Vries type, {\it viz.}
\begin{eqnarray}\label{coupledkdv}
\left\{
\begin{array}{lll}
u_t + u_{xxx} + P(u,v)_x = 0, \\
v_t + v_{xxx} + Q(u,v)_x = 0,
\end{array}
\right.
\end{eqnarray}
where $P$ and $Q$ are, say, homogeneous polynomials.  

Concerning singularity formation, it is an open question whether or not smooth, 
rapidly decaying,
real--valued solutions of (\ref{gkdv}) develop singularities in
finite time in the supercritical
case $p\ge 5$.   In the critical case $p = 4$, blowup in finite 
time has been established by Martel and Merle \cite{M}, \cite{MM} whilst 
for subcritical values $p = 1,2,3$, there is no singularity formation for 
data that lies at least in the Sobolev space $H^1(\R)$.  (However, 
 solutions corresponding to infinitely smooth initial values lying only in $L^2(\R)$ 
can develop singularities; see \cite{BS}.)    
Numerical simulations reported in \cite{BDKM} of solutions of
 (\ref{gkdv})  
initiated with an amplitude--modified solitary 
wave indicate blowup in finite time.  Such initial data has an 
analytic extension to a strip symmetric about the real axis in 
the complex plane.  The results of Bona, Grujic and Kalisch \cite{GrKa}, \cite{BoGr} and  \cite{BoGrKa} 
indicate that blowup at time $t$ has to be accompanied by the width of the 
strip of analyticity shrinking to zero at the same time.  This points to the 
prospect of a pair of complex conjugate singular points coliding 
at a spatial point on the real axis, thereby producing non--smooth behavior 
of the real--valued solution.  It was shown in \cite{BW2} that, in certain cases, 
curves of singularities  do merge together.  This 
happens  at the moment of interaction of the two solitons when the 
amplitudes are related in a particular way.  This result provides some 
indication 
that the blowup seen in the numerical simulations might occur because of
the coalescence of curves of complex singularities.   
Such ruminations seem to justify interest in the pole 
dynamics in the
context of (\ref{gkdv}), even in the case 
where the initial data is real--valued.

If one considers instead complex--valued solutions, it was shown in
\cite{BW1} (and see also \cite{BW3}) that, in the case of spatially periodic boundary conditions,
equation \re{gkdv} has solutions which blow up in finite time
for all integers $p \ge 1$.  Explicit examples of smooth, complex--valued solutions 
defined for $x, t \in \R$
of the Korteweg--de Vries equation \re{kdv} that blow 
up in finite time have been given in \cite{Bir,BW2,Li,WY}.
One outcome of the present paper is an explicit example of a 
blowing--up soution of the modified
Korteweg--de Vries equation \re{mkdv1}.

For the system (\ref{coupledkdv}) where $P$ and $Q$ are 
homogeneous quadratic polynomials, 
conditions on the coefficients are known that imply global well--posedness 
for real--valued initial data $(u_0,v_0)$ (see \cite{BCW}).  
And, the pole dynamics investigated  in \cite{BW2} for the Korteweg--de Vries 
equation (\ref{kdv}) itself, (KdV--equation henceforth) 
reveals that the choice
\begin{equation}\label{quadratic}
P(u,v) = u^2 - v^2 \quad \text{and} \quad Q(u,v) = 2uv
\end{equation}
leads to a system (\ref{coupledkdv})  possessing solutions that 
develop singularities in finite time.  (This is the system that obtains
if complex--valued solutions of the KdV--equation are broken up
into real and imaginary parts.)


The latter result is obtained by a careful study of the pole dynamics of the
explicit two--soliton solution $U = U(z,t) = U(x+iy,t)$ of the KdV--equation in the complex 
$z$--plane.  It transpires that as a function of time, most of the 
 singularities of this exact solution, which are all poles, move vertically in the $y$--direction 
in the complex plane as well as 
propagating horizontally in the $x$--direction.  As a consequence, by 
choosing $y_0$ appropriately, one can arrange that the function 
$u(x,t) = U(x + iy_0,t)$ is a complex--valued solution of the KdV--equation 
that, at $t=0$, is infinitely smooth and decays to zero exponentially rapidly 
as $x \to \pm\infty$, but which blows up for a positive value $t > 0$.  
If we write $u(x,t) = v(x,t) + iw(x,t)$, then the pair $(v,w)$ is a solution 
of (\ref{coupledkdv}) with the choice (\ref{quadratic}) that starts at 
$(v_0,w_0) = (v(\cdot,0),w(\cdot,0))$ 
smooth and rapidly decaying, but which forms a singularity in finite time.
It is worth noting that by an appropriate choice of the particular two--soliton 
solution, the initial data $(v_0,w_0)$ can be taken to be arbitrarily small in
any of the usual function spaces used in the analysis of such equations.

The motivation for the current paper was to see if the 
phenomenon just described, i.e. curves of singularities merging together, can
occur for the modified KdV--equation, and if this behavior
might provide additional insight into the possible ways a singularity
is produced in nonlinear, dispersve wave equaions.  Indeed, it turns out
that the same phenomenon does occur.   In certain cases, curves
of singularities do converge together at one point.  On the other
hand, the behavior of these curves is not so different from what
occurs in the  context of the KdV--equation.


\section{The two--soliton solutions}\label{formula}
Preliminary analysis of the two--soliton solutions of the
 mKdV--equation
\begin{equation}\label{mkdv}
u_t+u_{xxx}+6u^2u_x=0,\quad x\in\C, t\in\R,\end{equation}
are set forth here in preparation for the investigation of their pole
dynamics.
We begin with a standard transformation enabling one to
express solutions of \re{mkdv} in terms of solutions of another equation.
Let $u = v_x$ where $u$ is a solution of \re{mkdv}.  Then $v$ satisfies the equation
$$
\frac{d}{dx}\left(v_t+v_{xxx}+2v_x^3\right)=0.
$$
Assume $v$ and its derivatives vanish at infinity and search for solutions of 
the latter equation of the form $v = 2\arctan(g)$. 
A calculation shows that $v$ satisfies
\begin{equation}\label{eqv}
v_t+v_{xxx}+2v_x^3=0
\end{equation}
 if and only if
\begin{equation}\label{eqg}
(1+g^2)(g_t+g_{xxx})+6g_x(g_x^2-gg_{xx})=0.
\end{equation}
This yields a solution to \re{mkdv} having the form
\begin{equation}\label{mkdvgsol}
u(x,t)=2\big(\arctan g(x,t)\big)_x =  \frac{2 g_x(x,t)}{1 + g(x,t)^2}.
\end{equation}
It is important to note that equation \re{mkdv} and \re{eqg} are both
invariant under change of sign.   That is,   $u$ is a solution of \re{mkdv} if and
only if $-u$ is a solution and, likewise, $g$ is a solution of \re{eqg} if and
only if $-g$ is a solution.  Also, replacing $g$ by $1/g$ in
\re{mkdvgsol} has the effect of multiplying the solution $u$ by $-1$.
More precisely, if $u$ is given by \re{mkdvgsol}, then
\begin{equation}\label{mkdvgsolinv}
2\left(\arctan \frac{1}{g(x,t)}\right)_x = -2 \frac{g_x(x,t)}{1 + g(x,t)^2} = -u(x,t).
\end{equation}

The well--known soliton solution of (\ref{mkdv}) has a hyperbolic secant profile.  
In detail, for any amplitude value $k>0$, it is straightforward to check that
\begin{equation}\label{mkdvgsoliton}
g(x,t)=\exp\big(-k(x-x_0)+k^3t\big) = \exp\big(-k(x-x_0-k^2t)\big)
\end{equation}
is a solution of \re{eqg}. The corresponding solution $u$ of \re{mkdv} is the soliton solution with speed $k^2$ and is  given explicitly as
\begin{align}\label{solsol}
u(x,t)&=2\big(\arctan g(x,t)\big)_x = -2k\frac{e^{-k(x-x_0)+k^3t}}{1+e^{-2k(x-x_0)+2k^3t}} \\&= -k\sech\big(-k(x-x_0) +k^3t\big).
\end{align}
Note that the choice of sign in the exponential function $g$ in \re{mkdvgsoliton} produces
the {\it negative} soliton.  Replacing $g$ by either $-g$ or $1/g$ will produce
the positive soliton. 

The two--soliton solutions are a little more complicated.  The formulation presented here 
is based on that appearing in \cite{Bi}.  As just noted, there are both positive and negative soliton solutions of 
\re{mkdv}.  Consequently, there are two types of two--soliton solutions, namely 
interacting solitons of the same or opposite signs.

Suppose first that $0 < k_1 < k_2$.  Define the functions $f_j$ by
\begin{equation}\label{deffj}
f_j(x,t) = \exp(-k_jx+k_j^3t),\quad j=1,2.
\end{equation}
Of course, this definition omits two arbitary spatial translations;  these will be added later. 
Define two auxilliary functions, $g^+$ and $g^-$ by 
\begin{align}
\label{defg+}g^+(x,t) &= -\gamma\frac{f_1(x,t)+ f_2(x,t)}{1- f_1(x,t)f_2(x,t)},\\
\label{defg-}g^-(x,t) &= \gamma\frac{f_1(x,t)- f_2(x,t)}{1+ f_1(x,t)f_2(x,t)},
\end{align}
where
\begin{equation}\label{gamma}
\gamma=\frac{k_2+k_1}{k_2-k_1}>1.
\end {equation}

\begin{proposition}
The functions $g^+$ and $g^-$ defined in \re{defg+}-\re{defg-} are solutions to (\ref{eqg}).
\end{proposition}

\begin{proof}
It suffices to provide the proof for $g^-$.
Indeed, if one replaces $f_1$ by $-f_1$, then $g^-$ is transformed
into $g^+$, and all the calculations below remain valid in this case.

Temporarily, set $g=g^-$.
Notice that $f_{jx}=-k_jf_j$ and $f_{jt}=k_j^3f_j$ for $j=1,2$.  Thus, the quantities $g_t$, $g_x$, $g_{xx}$ and $g_{xxx}$ may be expressed in terms of $f_j$ and $k_j$.  First, differentiate with respect to time and come to the expression
$$g_t = \gamma\frac{k_1^3f_1-k_2^3f_2+k_1^3f_1f_2^2-k_2^3f_1^2f_2}{(1+f_1f_2)^2}.$$
Similarly, the derivative with respect to $x$ is
$$g_x = \gamma\frac{k_2f_2-k_1f_1-k_1f_1f_2^2+k_2f_1^2f_2}{(1+f_1f_2)^2}.$$
Differentiating the latter expression leads to   
$$g_{xx} = \frac{\gamma}{(1+f_1f_2)^3}\Big[k_1^2f_1-k_2^2f_2
+(k_1^2+4k_1k_2+k_2^2)(f_1f_2^2-f_1^2f_2)-k_1^2f_1^2f_2^3+k_2^2f_1^3f_2^2\Big].$$
Differentiating once more and simplifying gives
\begin{multline*}
g_{xxx} = \frac{\gamma}{(1+f_1f_2)^4}\Big [ k_2^3f_2-k_1^3f_1-(k_1^3+4k_2^3+6k_1^2k_2+12k_1k_2^2)(f_1f_2^2+f_1^3f_2^2) \\+ (4k_1^3+k_2^3+6k_1k_2^2+12k_1^2k_2)(f_1^2f_2+f_1^2f_2^3)-k_1^3f_1^3f_2^4+k_2^3f_1^4f_2^3)\Big].
\end{multline*}
It follows that 
$$g_t+g_{xxx} = \frac{6\gamma(k_1+k_2)^2}{(1+f_1f_2)^4} \big[ k_1f_1^2f_2-k_2f_1f_2^2+k_1f_1^2f_2^3-k_2f_1^3f_2^2 \big]$$
on the one hand, and 
\begin{align*}
g_x^2-gg_{xx} =  \frac{\gamma^2}{(1+f_1f_2)^4} \big[(k_1-k_2)^2(f_1f_2+f_1^3f_2^3)-  \\ 8k_1k_2f_1^2f_2^2+(k_1+k_2)^2(f_1^3f_2+f_1f_2^3) \big]
\end{align*}
on the other.  
It is now straightforward to check that $g$ satisfies equation (\ref{eqg}).

\end{proof}

The above proposition implies that
$$u^\pm(x,t)=2\big(\arctan g^\pm(x,t)\big)_x = \frac{2g^\pm_x(x,t)}{1+(g^\pm)^2(x,t)}$$
are solutions to \re{mkdv}.
A further calculations reveals that
\begin{equation}\label{uplus}
u^+=2\gamma \frac {G^+}{F^+}
\end{equation}
 and
 \begin{equation}\label{uminus}u^-= 2\gamma \frac {G^-}{F^-}
 \end{equation}
  where the new combinations 
\begin{align}
G^+ &=k_1f_1(1 + f_2^2) + k_2f_2(1 + f_1^2),\label{defG+}\\
G^- &= -k_1f_1(1 + f_2^2) + k_2f_2(1 + f_1^2),\label{defG-}\\
F^+ &= (1- f_1f_2)^2+\gamma^2(f_1+ f_2)^2,\label{defF+}\\
F^- &= (1+ f_1f_2)^2+\gamma^2(f_1- f_2)^2,\label{defF-}
\end{align}
have been introduced.  Note that the functions with a superscript ``$+$" are obtained from the functions
with a superscript ``$-$" simpy by replacing $f_1$ by $-f_1$.
If every occurence
of $k_1$ is replaced by $-k_1$  in formula \re{uminus} for $u^-$, then 
$f_1$ is replaced by $1/f_1$ and $\gamma$ is replaced by $1/\gamma$.
Simplifying the resulting expression gives exactly the
formula \re{uplus} for $u^+$.  In other words, the formulas for
$u^+$ and $u^-$ can be obtained from each other by replacing
every occurence of $k_1$ by $-k_1$.

The function $u^+$ is a two--soliton solution of \re{mkdv} having two positive
interacting solitons whist $u^-$ is a two--soliton solution with
 two interacting solitons of opposite sign, the faster one being the positive one.
One can see this graphically using MAPLE or Mathematica, for example.  Analytically, 
 there is an explicit relationship between the fomulas
for $u^\pm$ given by \re{uplus} and \re{uminus} and the
single soliton solutions of speeds $k_1^2$ and $k_2^2$.
As in \re{solsol}, let $u_j$
be the centered, positive, soliton solution 
\begin{equation}
u_j(x,t) = k_j\sech(-k_jx + k_j^3t) = \frac{2k_jf_j(x,t)}{1 + f_j(x,t)^2}.
\end{equation}
 of speed $k_j^2, j= 1,2.$ 
If the formulas for \re{uplus} and \re{uminus} are both divided
by $(1 + f_1^2)(1 + f_2^2)$, there obtains
\begin{equation}\label{sum}
u^+ = \gamma\frac{u_1 + u_2}{D^+}
\end{equation}
and
\begin{equation}\label{diff}
u^- = \gamma\frac{-u_1 + u_2}{D^-}
\end{equation}
where
\begin{align*}
D^+ &= \frac{(1- f_1f_2)^2+\gamma^2(f_1+ f_2)^2}{(1 + f_1^2)(1 + f_2^2)}, \\
D^- &= \frac{(1+ f_1f_2)^2+\gamma^2(f_1- f_2)^2}{(1 + f_1^2)(1 + f_2^2)}.
\end{align*}
Formulas \re{sum} and \re{diff} show in particular
that $u^+(x,t) > 0$ for all $x \in \R$ and that $u^-(x,t) > 0$ precisely for those $x \in \R$ and $t \in \R$
for which $u_2(x,t) > u_1(x,t)$.

  In Section \ref{largetime}, the 
asymptotic behavior of the singularities
of $u^\pm$ are examined for large positive and negative time.  We 
will see that they
separate into two groups, corresponding to the two single
solitons.  More remarks
on the shape of the two--soliton solutions during their interaction
are to be found in Section \ref{maple}.

Before ending this section, we return to the issue
of spatial shifts.
For $0<k_1<k_2$ and $x_1,x_2\in\R$, let
\begin{align*}
\tilde{f_1}(x,t)&=\exp(-k_1(x-x_1)+k_1^3t),\\
\tilde{f_2}(x,t)&=\exp(-k_2(x-x_2)+k_2^3t),
\end{align*}
and
\begin{align*}\widetilde{g^+}(x,t)=-\frac{\tilde{f_1}(x,t)+\tilde{f_2}(x,t)}{1-\gamma^{-2}\tilde{f_1}(x,t)\tilde{f_2}(x,t)},  \\
\widetilde{g^-}(x,t)=\frac{\tilde{f_1}(x,t)-\tilde{f_2}(x,t)}{1+\gamma^{-2}\tilde{f_1}(x,t)\tilde{f_2}(x,t)}.
\end{align*}
Define the {\it interaction time} $t_0$ and the {\it interaction center} $x_0$ 
of $\widetilde{g^\pm}(x,t)$ to be 
\begin{align}
\label{deft0}t_0 &= -\frac{x_2-x_1}{k_2^2-k_1^2}-\frac{1}{(k_2+k_1)k_1k_2}\log\gamma \quad {\rm and}\\
\label{defx0}x_0 &= \frac{k_2^2x_1-k_1^2x_2}{k_2^2-k_1^2}-\frac{k_1^2+k_1k_2+k_2^2}{(k_2+k_1)k_1k_2}\log\gamma,
\end{align}
respectively.

\begin{proposition}\label{prop-sym}  Let $t_0$ and $x_0$ be the interaction time and interaction
center for $\widetilde{g^+}(x,t)$ and $\widetilde{g^-}(x,t)$.   Then, for any $t\in\R$ and $x\in\C$, we have
\begin{equation}\label{eq-sym}\widetilde{g^\pm}(x,t) = g^\pm(x-x_0, t-t_0),\end{equation}
where $g^\pm$ are defined in \re{defg+}-\re{defg-}.
Moreover, the functions $\widetilde{u^\pm}(\cdot,t_0)=2(\arctan \widetilde{g^\pm}(\cdot, t_0))_x$ are symmetric about the point $x_0$ on both $\R$ and $\C$.
\end{proposition}
\begin{proof}
It suffices to find $(x_0,t_0)\in\R^2$ such that $\tilde{f_j}(x,t)=\gamma f_j(x-x_0, t-t_0)$ ($j=1,2$) for all $t\in\R$, $x\in\C$. Equivalently, (\ref{eq-sym}) will be satisfied if
$$
\left\{\begin{array}{ll}\gamma e^{k_1x_0-k_1^3t_0}=e^{k_1x_1}, \\ \gamma e^{k_2x_0-k_2^3t_0}=e^{k_2x_2}. \end{array}\right.
$$
Since $(x_0, t_0)$ given by (\ref{deft0})-(\ref{defx0}) is the solution for this system, the first assertion is proved.
To see the symmetry of the functions $\widetilde{u^\pm}(\cdot, t_0)$ about $x_0$, it is only 
necessary to deduce from (\ref{eq-sym}) that  $u^\pm(\cdot,0)=2(\arctan(g^\pm(\cdot,0))_x$ is 
an even function.  This is easily verified since $g^\pm(-x,0) = -g^\pm(x,0)$ for all $x\in\C$.
\end{proof}

Since (mKdV) is invariant under time-- and space--translation, it is concluded 
from Proposition \ref{prop-sym} that the functions 
$\widetilde{u^\pm}(x,t)=2(\arctan(\widetilde{g^\pm}(x,t))_x$ are also solutions to \re{mkdv}.
The point is that, for a general two--soliton solution, the time and
place of the interaction are given by $t_0$ and $x_0$ in \re{deft0} and \re{defx0}, 
respectively.  In particular, the solutions ${u^\pm}$ are already
normalized so that the interaction time is $t = 0$, and the center
of the interaction is at $x = 0$.  It is interesting to note that
the values of $t_0$ and $x_0$ in the above
proposition are precisely the same as for the two--soliton solution
of the Korteweg deVries equation, as given in Theorem 1 of \cite{BW2}.

\section{Singularities of the two--soliton solutions in the complex plane}
As mentioned in the introduction, the two--soliton solutions $u(x,t)$ of \re{mkdv} are viewed as 
meromorphic functions in the complex variable $x$.  
We aim to determine how the dynamics of the singularities  in $\C$
 reflect the behavior of $u(x,t)$ when $x$ is restricted 
to the real axis.

Consider first 
the (one)--soliton solution, given by \re{solsol}.
It is immediate that the singularities of these solutions are 
precisely 
\begin{equation}\label{solsing}
x = x_0 + k^2t + \frac{m\pi i}{2k}
\end{equation}
where $m$ runs through the odd integers.  Moreover, these singularities of $u$ are all simple
poles.  Thus, the speed  of the soliton is exactly
 the speed of its poles in the complex plane, while the position of the maximum point of the soliton
at time $t$  is  the real part of the position of these poles.  The imaginary part of
the singularity remains constant in time of course.  

The poles of the two--solitons solutions \re{uplus}
and \re{uminus}, correspond to the zeros of the functions
$F^\pm$ and $G^\pm$ defined in \re{defG+}, \re{defG-}, \re{defF+}
and \re{defF-}.  It will turn out that the singularities of $u^\pm$
are all simple poles, just as for the one--solitons.  In all but a specfic
 class of exceptional cases, these poles 
correspond to simple zeros of $F^\pm$.

Recall that $k_1$ and $k_2$ are called {\it commensurable} if there
 exist positive integers $p_1$ and $p_2$ such that
\begin{equation}
\label{defCom}\frac{k_2}{k_1}=\frac{p_2}{p_1}.  
\end{equation}
Without loss of generality, we may take it that $\gcd(p_1, p_2)=1$.  In what follows, we will
always assume that the integers $p_1$ and $p_2$ appearing in \re{defCom} are without common
prime factors. In this case, $F^\pm(x,t)$ and $G^\pm(x,t)$ (and thus $u^\pm(x,t)$) are periodic in $x$ with minimal imaginary period $2\pi\lambda i$ where
\begin{equation}\label{defLam}\lambda = \frac{p_1}{k_1} = \frac{p_2}{k_2}.\end{equation}

\begin{lemma}\label{lem-multzeros}
For any $t\in\R$, the zeros of $F^\pm(\cdot, t)$ are simple, except for the following special case. If $k_1$ and $k_2$ are commensurable, and if $p_1,p_2\in\N$ and $\lambda>0$ are given by \re{defCom} and \re{defLam}, with $p_1\pm p_2\in 4\N$,
then there is a fourth order zero of $F^\pm(\cdot, 0)$ at $x=(\frac 12+q)\lambda\pi i$ for all $q\in\Z$.
\end{lemma}

\begin{proof}
From \re{defF-}, we infer
$$F^-_x = 2\big[(1+f_1f_2)(f_{1x}f_2+f_1f_{2x})+\gamma^2(f_1-f_2)(f_{1x}-f_{2x})\big]. $$
Noticing that $f_{jx} = -k_jf_j$, ($j=1,2$) and setting $X=f_1$ and $Y=f_2$, 
it follows that if $x$ is a zero of $F^-(\cdot,t)$ of order greater than or equal to 2, then
\begin{equation}\label{syst-XY}
\left\{\begin{array}{ll}(1+XY)^2+\gamma^2(X-Y)^2=0 \quad \text{and}\\ -(k_1+k_2)(1+XY)XY+\gamma^2(X-Y)(k_2Y-k_1X)=0.\end{array}\right.
\end{equation}
Let $(X,Y)$ be a solution of this system. Then, from the first equation we deduce that
\begin{equation}\label{eqXY1}
1+XY = i\eps\gamma(X-Y),
\end{equation}
where $\epsilon$ is either 1 or -1.  
Inserting this into the second equation, it is found that
$$i\eps(k_2-k_1)XY = k_2Y-k_1X,$$
provided $X\neq Y$. Extract the product $XY$ and inject it in (\ref{eqXY1}) to obtain the linear relation
$$Y = \frac{k_2}{k_1}X + i\eps\frac{k_2-k_1}{k_1}.$$
Then formla (\ref{eqXY1}) implies that $X$ must satisfy
$$1+X\left(\frac{k_2}{k_1}X+i\eps\frac{k_2-k_1}{k_1}\right) = i\eps\gamma\left(\left(1-\frac{k_2}{k_1}\right)X-i\eps\frac{k_2-k_1}{k_1}\right),$$
which simplifies to
$$X^2+2i\eps X-1 = 0.$$

It follows that $X=Y=\pm i$. Thus, it must be that $e^{-k_1x+k_1^3t}=e^{-k_2x+k_2^3t}=\pm i$, from which we deduce that $-k_1x+k_1^3t$ and $-k_2x+k_2^3t$ are both purely imaginary.
Therefore, $t=0$ and $x$ is purely imaginary. Additionally, we have
$$\left\{\begin{array}{ll}k_1x=(\frac 12+q_1)\pi i, \\ k_2x=(\frac 12+q_2)\pi i, \end{array}\right. \quad \text{for some  } 
q_1,q_2\in\Z.$$
Thus $k_1$ and $k_2$ are commensurable and in fact, $k_2/k_1=r_2/r_1$ with $r_j=1+2q_j$. Noticing that $q_2-q_1\in2\Z$ (because $k_1x-k_2x\in2i\pi\Z$), it transpires that $r_2-r_1=4k$ for some $k\in\Z$. It follows that $d=\gcd(r_1,r_2)$ divides $4k$ and therefore $d$ divides $k$ since $d$ is necessarily an odd integer. Define $p_1$ and $p_2$ by $k_2/k_1=p_2/p_1$ with $\gcd(p_1,p_2)=1$. With this 
definition, $p_2-p_1=(r_2-r_1)/d = 4k/d\in4\N$.
 There are precisely two such complex numbers in the fundamental strip $S=\{x\in\C, -\lambda\pi<\Im x< \lambda\pi\}$, and they are $x=\pm\frac\lambda 2\pi i$.

It is straightforward to ascertain that for such values of $x$, 
  $F^-_{xx}(x,0)=F^-_{xxx}(x,0)=0$ and $F^-_{xxxx}(x,0)\neq 0$. 
The same arguments hold for $F^+$.
\end{proof}

\begin{lemma}\label{lem-Gzeros}
For any $t\in\R$, the zeros of $F^-(\cdot, t)$ and $G^-(\cdot, t)$ are distinct, except the following special case. If $k_1$ and $k_2$ are commensurable, and if $p_1,p_2\in\N$ and $\lambda>0$ are given by \re{defCom} and \re{defLam}, with $p_1\pm p_2\in 4\N$,
then there is a third order zero of $F^-(\cdot, 0)$ at $x=(\frac 12+q)\lambda\pi i$ for all $q\in\Z$.  The same statement is true with  $F^-$ and $G^-$ replaced by $F^+$ and $G^+$.
\end{lemma}

\begin{proof}
In the notation of the proof of Lemma \ref{lem-multzeros}, the result follows if the system
\begin{equation}\label{syst-XY2}
\left\{\begin{array}{ll}(1+XY)^2+\gamma^2(X-Y)^2=0,\\ -k_1X(1+Y^2)+k_2Y(1+X^2)=0,\end{array}\right.
\end{equation}
admits as its only solution $X=Y=\pm i$. As before, if $(X,Y)$ is a solution, there exists $\eps\in\{-1,1\}$ for which
$$1+XY = i\eps\gamma(X-Y).$$
Extracting from this the variable $Y$ and computing $1+Y^2$ leads to
$$Y=\frac{i\eps\gamma X-1}{X+i\eps\gamma}$$
and so
$$1+Y^2=(1-\gamma^2)\frac{1+X^2}{(X+i\eps\gamma)^2}.$$
Now insert this into the second equation of the system \re{syst-XY2} and, assuming by 
contradiction that $1+X^2$ is not zero, simplify the outcome.  It follows that
$$\frac{k_1}{k_2}(1-\gamma^2)X = (i\eps\gamma X-1)(X+i\eps\gamma).$$
After further simplifications, this becomes
$$X^2+2i\eps X-1=0,$$
and the claim follows.
\end{proof}
\begin{definition}\label{exc}  The values of $k_1$ and $k_2$ for which
 $F^\pm(\cdot,0)$ has multiple zeros are collectively referred to as the {\rm exceptional case}. This occurs when $k_1$ and $k_2$ are commensurable,
$p_1$ and $p_2$ are odd, and
$$p_2-p_1\in 4\N$$
when considering $F^-$, and
$$p_2+p_1\in 4\N$$
when considering $F^+$.
\end{definition}

Lemmas \ref{lem-multzeros} and \ref{lem-Gzeros} show that all
singularities of $u^\pm(\cdot,t)$ are simple poles.  In all but the exceptional
case, they correspond to simple zeros of $F^\pm$.  In the exceptional case,
they correspond to a third order zero of $g^\pm$ coinciding with a
fourth order zero of $F^\pm$.  It follows from Rouch\'e's Theorem that
in the exceptional case, four simple poles converge, as $t \to 0$, to the
simple pole at $(\frac 12+q)\lambda\pi i$ for all $q\in\Z$.
Moreover, by the residue theorem, the sum of the residues of
$u^\pm(\cdot,t)$ at these four poles converge, as $t \to 0$, to the
residue of the simple pole at $(\frac 12+q)\lambda\pi i$ for all $q\in\Z$.

\begin{proposition}\label{curves}
In all but the exceptional case, the poles of $u^\pm(\cdot,t)$, i.e. the zeros of
$F^\pm(\cdot,t)$, are described by analytic curves $x : \R \to \C$.  In the exceptional case, these curves are defined and analytic separately
for $t < 0$ and $t > 0$.
\end{proposition}

\begin{proof}
Consider the case of $u^+$ and $F^+$.
Since all the zeros of $F^+(\cdot,t)$ are simple, the implicit--function theorem 
shows that, for a fixed time $t_0$, a zero of $F^+(\cdot,t_0)$ can be locally 
and uniquely continued as an analytic curve $x(t)$ such that $F^+(x(t),t)  = 0$.
Such a curve $x(t)$ can be continued as long as it remains in a bounded
region of $\C$.  Thus, we need to show that $|x(t)|$ must stay bounded
as long as $t$ remains in a bounded interval of $\R$.  First, it follows
from Proposition \ref{bddimpt} in Section \ref{vert} below that the imaginary part
of $x(t)$ must remain bounded.  Furthermore, if $\Re x(t) \to \infty$ in finite time, then
the equation $F^+(x(t),t)  = 0$ implies $1 = 0$.  If $\Re x(t) \to -\infty$
in finite time, then the equation $\frac{F^+(x(t),t)}{1 - f_1(x(t),t)^2f_2(x(t),t)^2}  = 0$ also implies that $1 = 0$.

A similar argument works for $u^-$ and $F^-$.

\end{proof}

This section is closed with a specific example of an element in the exceptional case, namely $k_1=1$, $k_2=5$ (calculations done with MAPLE).
In this case,  $\lambda=1$ and $F^-$ and $G^-$ may be rewritten as
\begin{align*}
G^-(x,t) &= -e^{251t}y^{11}+5e^{127t}y^7+5e^{125t}y^5-e^ty \text{  and  } \\
F^-(x,t) &= e^{252t}y^{12}+\frac 94e^{250t}y^{10}-\frac 52e^{126t}y^6+\frac 94e^{2t}y^2+1,
\end{align*}
where $y=e^{-x}$. Note that at time $t=0$, $F^-$ and $G^-$ are symmetric polynomials in $y$,
 and that $i$ and $-i$ are third order zeros of $G^-$, and fourth order zeros of $F^-$.  All the 
 other zeros are simple.
We also have the decomposition
$$u^-(x,0) = 3\frac {G^-(x,0)}{F^-(x,0)} = \frac{-2}{y-i}-\frac{2}{y+i}+\frac{y-4i}{2y^2-iy-2}+\frac{y+4i}{2y^2+iy-2}.$$

\section{Large--time asymptotic behavior of the singularities}\label{largetime}
In this section we show that the poles of $u^\pm$,
that is the zeros of $F^\pm$, separate out into two groups,
as $t \to \pm\infty$, behaving asymptotically as poles of
single solitons, of speeds $k_1^2$ and $k_2^2$, respectively.
This, of course, reflects the fact that $u^\pm$ are in fact
``two--soliton" solutions of \re{mkdv}.

The situation is nearly identical to that obtaining for 
the two--soliton solution for the Korteweg--de Vries
equation (compare the following theorem with Theorem 2 in \cite{BW2}).
In particular, the backward shift of the slower wave, well--known
in the case of the KdV--solitons, is also present for the modified KdV--equation.

\begin{theorem}\label{asympt}
 The asymptotic behaviors as $t \to \pm\infty$ of  the curves $x(t)$ of zeroes 
 of $F^\pm(x, t)$ whose existence was determined in Proposition \ref{curves} are completely described as follows.
\begin{enumerate}
 \item For every odd integer $m\in\Z$, there exists a unique curve $x_{m, s^-}(t)$ of zeros of $F^\pm(\cdot, t)$ such that
 $$x_{m, s^-}(t) = k_1^2t+\frac 1{k_1}\log\gamma+\frac{m\pi i}{2k_1}+o(1),$$
 as $t\to-\infty$.
 \item For every odd integer $m\in\Z$, there exists a unique curve $x_{m, s^+}(t)$ of zeros of $F^\pm(\cdot, t)$ such that
 $$x_{m, s^+}(t) = k_1^2t-\frac 1{k_1}\log\gamma-\frac{m\pi i}{2k_1}+o(1),$$
 as $t\to\infty$.
 \item For every odd integer $n\in\Z$, there exists a unique curve $x_{n, f^-}(t)$ of zeros of $F^\pm(\cdot, t)$ such that
 $$x_{n, f^-}(t) = k_2^2t-\frac 1{k_2}\log\gamma+\frac{n\pi i}{2k_2}+o(1),$$
 as $t\to-\infty$.

\item For every odd integer $n\in\Z$, there exists a unique curve $x_{n, f^+}(t)$ of zeros of $F^\pm(\cdot, t)$ such that
 $$x_{n, f^+}(t) = k_2^2t+\frac 1{k_2}\log\gamma-\frac{n\pi i}{2k_2}+o(1),$$
 as $t\to\infty$.\end{enumerate}

\end{theorem}

To prove Theorem \ref{asympt}, it is necessary to
study the zeros of $F^\pm(\cdot,t)$ with respect to
frames of reference which move at the speed of each
constituant soliton.  As in Section 4 of \cite{BW2}, we set
\begin{align}
\label{MovFraz}z &= x-k_1^2t,\\
\label{MovFraw}w &= x-k_2^2t,\\
\label{MovFrar}r &= \exp(k_2(k_2^2-k_1^2)t),\\
\label{MovFras}s &= \exp(k_1(k_2^2-k_1^2)t).
\end{align}
If
\begin{equation}\label{H}
H^\pm(z,r)=(1\mp re^{-(k_1+k_2)z})^2+\gamma^2(e^{-k_1z}\pm re^{-k_2z})^2
\end{equation}
and
\begin{equation}
 I^\pm(w,s)= (s\mp e^{-(k_1+k_2)w})^2+\gamma^2(e^{-k_1w}\pm se^{-k_2w})^2,
\end{equation}
then
$$F^\pm(x,t) = H^\pm(z,r) = s^{-2}I^\pm(w,s).$$

Solutions $z(r)$ of $H^\pm(\cdot,r)$ which remain bounded
in $\C$ as $r \to 0$ and as $r \to \infty$ correspond to zeros $x(t)$
of $F^\pm(\cdot,t)$ asymptotically traveling at speed $k_1^2$
as $t \to \pm\infty$. Likewise,  solutions $w(s)$ of $I^\pm(\cdot,s)$ which remain bounded
in $\C$ as $s \to 0$ and as $s \to \infty$ correspond to zeros $x(t)$
of $F^\pm(\cdot,t)$ asymptotically traveling at speed $k_2^2$
as $t \to \pm\infty$.  The following is true of the curves $z(r)$.

\begin{proposition}
 For every odd integer $m\in\Z$, there exists a smooth curve $z_m^\pm(r)$ defined in some interval of $r\ge 0$, such that $H^\pm(z_m^\pm(r), r)=0$ and
\begin{equation}
 z_m^\pm(r) = \frac 1{k_1}\log\gamma +\frac{m\pi i}{2k_1}\mp(-1)^{\frac{m-1}{2}}\frac{4k_2}{k_2^2-k_1^2}\gamma^{-k_2/k_1}e^{-\frac{k_2}{2k_1}m\pi i}r+o(r)
\end{equation}
as $r\to 0^+$. For every odd integer $n\in\Z$, there exists a smooth curve $w_n^\pm(s)$ defined in some interval of $s\ge 0$, such that $I^\pm(w_n^\pm(s), s)=0$ and
\begin{equation}
 w_n^\pm(s) = -\frac 1{k_2}\log\gamma +\frac{n\pi i}{2k_2}\mp(-1)^{\frac{n-1}{2}}\frac{4k_1}{k_2^2-k_1^2}\gamma^{-k_1/k_2}e^{\frac{k_1}{2k_2}n\pi i}s+o(s)
\end{equation}
as $s\to 0^+$.
\end{proposition}

\begin{proof}
The relation $H^\pm(z,0)=0$ is satisfied if and only if there exists an odd integer $m\in\Z$ such that
$$z = \frac 1{k_1}\log\gamma + \frac{m\pi i}{2k_1},$$
and similarly, $I^\pm(w,0)=0$ is equivalent to
$$w=-\frac{1}{k_2}\log \gamma+\frac{n\pi i}{2k_2}$$
for some odd integer $n\in\Z$. Applying the implicit function theorem, there exist smooth curves $z_m^\pm(r)$ and $w_n^\pm(s)$ defined in a neighborhood of $z_0=1/k_1\log\gamma+m\pi i/2k_1$ and $w_0=-1/k_2\log\gamma+n\pi i/2k_2$, respectively, such that $H^\pm(z_m^\pm(r),r)=0$, $z_m^\pm(0)=z_0$ and $I(w_n^\pm(s),s)=0$, $w_n^\pm(0)=w_0$. It remains to calculate $(z_m^\pm)'(r)$ and $(w_n^\pm)'(r)$. Differentiating the equation $H^\pm(z_m^\pm(r),r)=0$ with respect to $r$ yields
\begin{multline*}
 \! \! \! \! \!  \big[-(k_1+k_2)(1\mp re^{-(k_1+k_2)z})re^{-(k_1+k_2)z}+\gamma^2(e^{-k_1z}\pm re^{-k_2z})(k_2re^{-k_2z}\pm k_1e^{-k_1z})\big]z'\\
 = -(1\mp re^{-(k_1+k_2)z})e^{-(k_1+k_2)z}+\gamma^2(e^{-k_1z}\pm re^{-k_2z})e^{-k_2z}
\end{multline*}
where $z=z_m^\pm(r)$. Taking $r=0$ gives
$$\pm k_1\gamma^2 e^{-2k_1z_m(0)}(z_m^\pm)'(0) = (\gamma^2-1)e^{-(k_1+k_2)z_m^\pm(0)},$$
and using that $$e^{-(k_1+k_2)z_m^\pm(0)} = \frac 1\gamma (-1)^{\frac{m-1}2}\gamma^{-k_2/k_1}e^{-\frac{k_2}{2k_1}m\pi i} $$ leads to 
$$(z_m^\pm)'(0)=\mp(-1)^{\frac{m-1}{2}}\frac{4k_2}{k_2^2-k_1^2}\gamma^{-k_2/k_1}e^{-\frac{k_2}{2k_1}m\pi i}.$$

The derivative $(w_n^\pm)'(0)$ is similarly calculated. Differentiating the equation $I^\pm(w_n^\pm(s),s)=0$ with respect to $s$ leads to
\begin{multline*}
 \! \! \!  [\pm(k_1+k_2)(s\mp e^{-(k_1+k_2)w})e^{-(k_1+k_2)w}+\gamma^2(e^{-k_1w}\pm se^{-k_2w})(\mp k_2se^{-k_2w}-k_1e^{-k_1s})]w'\\
 = -(s\mp e^{-(k_1+k_2)w})\mp\gamma^2(e^{-k_1w}\pm se^{-k_2w})e^{-k_2w}.
\end{multline*}
At $s=0$, there obtains
$$[-(k_1+k_2)e^{-2(k_1+k_2)w_n^\pm(0)}-k_1\gamma^2e^{-2k_1w_n^\pm(0)}](w_n^\pm)'(0) = \mp(\gamma^2-1)e^{-(k_1+k_2)w_n^\pm(0)}$$
and since
$$e^{-(k_1+k_2)w_n^\pm(0)} = -\frac 1\gamma (-1)^{\frac{n-1}2}\gamma^{-k_1/k_2}e^{-\frac{k_1}{2k_2}n\pi i}, $$
we conclude
$$(w_n^\pm)'(s) = \mp(-1)^{\frac{n-1}{2}}\frac{4k_1}{k_2^2-k_1^2}\gamma^{-k_1/k_2}e^{\frac{k_1}{2k_2}n\pi i}.$$
\end{proof}

A proof of Theorem \ref{asympt} is now readily available.
The last proposition immediately gives all the curves
of zeros of $F^\pm(\cdot,t)$ claimed in Theorem \ref{asympt}.
The only remaining issue is to prove that there are no
additional zeros of $F^\pm(\cdot,t)$.  If $k_1$ and $k_2$ are
commensurable, this is straightforward.  For each $t \in \R$, $F^\pm(\cdot,t)$
is a polynomial in $e^{-x\lambda}$ of degree $2(p_1 + p_2)$,
and thus must have precisely $2(p_1 + p_2)$ zeros in the complex
plane.
The zeros of $F^\pm(\cdot,t)$ described in Theorem \ref{asympt}
account for all of them, for large $|t|$, and so no other zeros can exist.

If $k_1$ and $k_2$ are not commensurable, the desired result
can be obtained by approximating $k_1$ and $k_2$ by sequences
 $\{k_1^\nu\}_{\nu = 1}^\infty$ and $\{k_2^\nu\}_{\nu = 1}^\infty$, which are commensurable.

\section{The nature of the singularity in the exceptional case}

In this section a more detailed analysis is undertaken of the singularity
of $u^\pm$ in the exceptional case (see Definition \ref{exc}).
As described just after this definition, the singularity
at $(\frac 12+q)\lambda\pi i$, for any  $q\in\Z$, corresponds
to a fourth order zero of $F^\pm(\cdot,0)$ and is approached
by four simple zeros of $F^\pm$ as $t \to 0$.  The goal is to understand the behavior of a smooth curve of zeros of $F(\cdot,t)$ in a neighborhood of such a fourth order zero.

We claim it suffices to analyse the singularity at $\frac\lambda 2\pi i$,
i.e. the case $q = 0$. Indeed, since $F^\pm(\cdot,t)$ is $2\pi\lambda i$ periodic, it 
is enough to consider $q = -1, 0$.  Furthermore,
$F^\pm(\overline x,t) = \overline{F^\pm(x,t)}$, and so
$F^\pm(x(t),t) = 0$ if and only if $F^\pm(\overline x(t),t) = 0$.
Hence,  only the case $q = 0$ need be examined.

Remark that
 $$
 F^\pm(x,t)=e^{-2(k_1+k_2)x+2(k_1^3+k_2^3)t}F^\pm(-x,-t).
 $$
 From this, it is deduced
that $F^\pm(x(t),t) = 0$ if and only if $F^\pm(-\overline x(t),-t) = 0$.
In other words, if $x(t)$ is a curve of zeros approaching $\frac\lambda 2\pi i$
as $t \nearrow 0$, then $-\overline x(t)$ is a curve
of zeros approaching $\frac\lambda 2\pi i$
as $t \searrow 0$.

\begin{theorem}\label{excsing}
Suppose we are in the exceptional case wherein $k_1$ and $k_2$ are commensurable and the odd integers $p_1, p_2\in\N$ and $\lambda>0$ are as in \re{defCom} and \re{defLam}.  In the case of $F^+$, it is 
assumed that  $p_1 + p_2\in 4\N$ and in the case of $F^-$, it is presumed that $p_1 - p_2\in 4\N$.
Let $x(t)$ be a smooth curve, defined for $t$ close to 0, such that $F^\pm(x(t),t)=0$ and $x(t)\neq \frac\lambda 2\pi i$ but $x(t)\to \frac\lambda 2\pi i$ as $t\to 0$.
Then, either
\begin{equation}\label{cubebeh}
\lim_{t\to 0}\frac{(x(t)-\frac\lambda 2\pi i)^3}{t} = -12,
\end{equation}
or
\begin{equation}\label{linbeh}
\lim_{t\to 0}\frac{x(t)-\frac\lambda 2\pi i}{t} = k_1^2+k_2^2.
\end{equation}
\end{theorem}

\begin{remark}
The similarity of the result in Theorem \ref{excsing} with the ``exceptional case" for the two--soliton
solution of sthe KdV--equation is quite striking.  The behavior of the curves
described in \ref{cubebeh} is exactly the same as given
by Proposition 4.9 in \cite{BW2}.  However, the exceptional case
for the  KdV--equation corresponds to $p_1$ being an odd integer, and $p_2$ being
an even integer.  For the mKdV--equation, the ``exceptional pole" is located at
$\frac{\lambda}{2} \pi i$, rather than at $ \lambda \pi i$ as it is for the KdV--equation. 

Also, how does one explain the behavior
described by \re{linbeh}?  In the case of the two--soliton solution of the
 KdV--equation, if $p_1$ is even and $p_2$ is odd, there is a horizontally
moving pole approaching the singularity at $\pi\lambda i$ exactly
as described by \re{linbeh}.  This calculation was not carried out in \cite{BW2},
but can be obtained by a simple modification of the proof of
Proposition 4.9 in \cite{BW2}.  Thus, it appears that the
exceptional case for the two--soliton solution of the  mKdV--equation \re{mkdv} includes
the behavior of curves of singularities from two different cases of the
two--soliton solution of the KdV--equation, namely the ``exceptional case", where
$p_1$ is odd and $p_2$ is even, as well as the case where
$p_1$ is even and $p_2$ is odd.  It is precisely in these two cases that there
is a pole located at the same place $\pi\lambda i$ (at time 0).
\end{remark}

\begin{proof}
The proof is provided for $F^-$.  Similar arguments apply for $F^+$.
Let $z(r)=x(t)-k_1^2t$ where $r$ is defined in \re{MovFrar}. It follows that
 $z(r)$ is a smooth curve such that $H^-(z(r),r)=0$, where $H$ is as in
\re{H}.  Notice that $z(r)\to \frac\lambda 2\pi i$ as $r\to 1$. Now, 
\begin{equation}\label{Hzero}
H^-(z,r)=0\Leftrightarrow 1+re^{-(k_1+k_2)z} = \pm i\gamma (e^{-k_1z}-re^{-k_2z}).
\end{equation}
If the minus sign is chosen on the right--hand side of \re{Hzero}, we come to
$$r = -\frac{i\gamma e^{-k_1z}+1}{e^{-(k_1+k_2)z}-i\gamma e^{-k_2z}}.$$
Differentiating this relation with respect to $r$, it follows that
\begin{align*}
1 = \frac{1}{ (e^{-(k_1+k_2)z}-i\gamma e^{-k_2z})^2 } \Big[ i\gamma k_1e^{-k_1z}\big(e^{-(k_1+k_2)z}-i\gamma e^{-k_2z}\big) \\  +\big(i\gamma e^{-k_1z}+1)(-(k_1+k_2)e^{-(k_1+k_2)z}+i\gamma k_2e^{-k_2z}) \Big] z'(r).
\end{align*}
This may be rewritten as
$$1 = \frac{i\gamma k_2e^{-k_2z} \big(1+i e^{-k_1z}\big)^2} {\big(e^{-(k_1+k_2)z}-i\gamma e^{-k_2z} \big)^2}z',$$
or what is the same,
\begin{equation}\label{zprime}
z' = -\frac{i e^{k_2z}  \big(e^{-(k_1+k_2)z}-i\gamma e^{-k_2z}\big)^2}{k_2\gamma \big(1+i e^{-k_1z}\big)^2}.
\end{equation}

Since $p_1$ and $p_2$ are both odd with $p_2 - p_1 \in 4\N$ and
$z(r) \to  \frac\lambda 2\pi i$ as $r \to 1$, it must be that
$e^{-k_1z} \to e^{-p_1\pi i/2}$ and $e^{-k_2z} \to e^{-p_2\pi i/2}$
as $r \to 1$.  If $p_1$ and $p_2$ are both in $4\N+ 1$, then
$e^{-k_1z}$ and $e^{-k_1z}$ both converge to $-i$ as $r \to 1$, while
if $p_1$ and $p_2$ are both in $4\N+ 3$, then
$e^{-k_1z}$ and $e^{-k_1z}$ both converge to $i$ as $r \to 1$.
We suppose first that $p_1$ and $p_2$ are both in $4\N+ 1$.
In this case, it follows from \re{zprime} that
\begin{equation}\label{zprime1}
z'(r) \to \frac{(\gamma+1)^2}{4k_2\gamma}=\frac{k_2}{k_2^2-k_1^2}
\end{equation}
as $r\to 1$, and thus
$$\lim_{r\to 1}\frac{z(r)-\frac\lambda 2\pi i}{r-1} = \frac{k_2}{k_2^2-k_1^2}.$$
Turning back to $x(t)$, and using the fact that $(r-1)/t\to k_2(k_2^2-k_1^2)$ as $t\to 0$, it follows that
$$\lim_{t\to 0}\frac{x(t)-k_1^2t-\frac\lambda 2\pi i}{t} = k_2^2$$
from which it is concluded that
$$\lim_{t\to 0}\frac{x(t)-\frac\lambda 2\pi i}{t} = k_1^2+k_2^2.$$

If instead the plus sign is chosen on the right--hand side of \re{Hzero}, then 
it is immediately inferred  that
$$r = \frac{i\gamma e^{-k_1z}-1}{e^{-(k_1+k_2)z}+i\gamma e^{-k_2z}}.$$
Differentiating this equation with respect to $r$, we get
\begin{align*}
1 = \frac{1}{ (e^{-(k_1+k_2)z}+i\gamma e^{-k_2z})^2  } \Big[ -i\gamma k_1e^{-k_1z}\big(e^{-(k_1+k_2)z}+i\gamma e^{-k_2z}\big) \\ -\big(i\gamma e^{-k_1z}-1\big)\big(-(k_1+k_2)e^{-(k_1+k_2)z}-i\gamma k_2e^{-k_2z}\big) \Big] z',
\end{align*}
which may be simplified to 
$$1 = -\frac{i\gamma k_2e^{-k_2z}(1-i e^{-k_1z})^2}{(e^{-(k_1+k_2)z}+i\gamma e^{-k_2z})^2}z'.$$
From this, it is inferred that
\begin{equation}\label{zprime2}
(1-i e^{-k_1z})^2z' = \frac{i e^{k_2z}}{k_2\gamma}(e^{-(k_1+k_2)z}+i\gamma e^{-k_2z})^2.
\end{equation}
Since $p_1$ and $p_2$ are both in $4\N+ 1$, it is concluded that
\begin{equation}\label{zprime3}
(1-i e^{-k_1z})^2z' \to -\frac{(\gamma-1)^2}{k_2\gamma}=\frac{-4k_1^2}{k_2(k_2^2-k_1^2)}
\end{equation}
as $r\to 1$. A consequence of this is that
$$\frac {d}{dr}(1-i e^{-k_1z})^3 = 3(1-i e^{-k_1z})^2(i k_1e^{-k_1z})z'\to \frac{-12k_1^3}{k_2(k_2^2-k_1^2)}.$$
L'Hopital's rule comes to the rescue and it is found that
$$\lim_{r\to 1}\frac{(1-i e^{-k_1z})^3}{r-1}= \frac{-12k_1^3}{k_2(k_2^2-k_1^2)}.$$
Since
$$\lim_{z\to \frac\lambda 2\pi i}\frac{1-i e^{-k_1z}}{z-\frac\lambda 2\pi i}=-i\lim_{z\to \frac\lambda 2\pi i}\frac{e^{-k_1z}-e^{- k_1\frac\lambda 2\pi i}}{z-\frac\lambda 2\pi i} = -i k_1e^{- k_1\frac\lambda 2\pi i} = k_1,$$
we must have
$$\lim_{r\to 1}\frac{(z(r)-\frac\lambda 2\pi i)^3}{r-1}=\frac{-12}{k_2(k_2^2-k_1^2)}.$$
Reverting to the original variable  $x(t)$  and using again that $(r-1)/t\to k_2(k_2^2-k_1^2)$ as $t\to 0$, it follows that
$$\lim_{t\to 0}\frac{(x(t)-k_1^2t-\frac\lambda 2\pi i)^3}{t} = -12$$
and thus
$$\lim_{t\to 0}\frac{(x(t)-\frac\lambda 2\pi i)^3}{t} = -12.$$

It remains to treat the situation where $p_1,p_2\in 4\N+3$, in which case $e^{-k_jz(r)}\to i$ as $r\to 1$. 
In fact, this case is ``dual" to the one just treated and the same calculations lead to the result. 
To see this, assume first that the minus sign obtains on the right--hand side of \re{Hzero}. Then \re{zprime} 
implies that
$$(1+i e^{-k_1z})^2z' = -\frac{i e^{k_2z}}{k_2\gamma}(e^{-(k_1+k_2)z}-i\gamma e^{-k_2z})^2\to -\frac{(\gamma-1)^2}{k_2\gamma}=\frac{-4k_1^2}{k_2(k_2^2-k_1^2)}$$
as $r\to 1$.
On the other hand, it is straightforward that
$$\lim_{z\to \frac\lambda 2\pi i}\frac{1+i e^{-k_1z}}{z-\frac\lambda 2\pi i}= k_1.$$
Following the same line of development as pursued for the positive sign in the previous  situation 
where both $p_1$ and $p_2$ lie in $4\N + 1$ leads to
$$\lim_{t\to 0}\frac{(x(t)-\frac\lambda 2\pi i)^3}{t} = -12.$$
Now assume we have the positive sign on the right--hand side of \re{Hzero}. From \re{zprime2} we deduce
$$z' = \frac{i e^{k_2z}}{k_2\gamma}\frac{(e^{-(k_1+k_2)z}+i\gamma e^{-k_2z})^2}{(1-i e^{-k_1z})^2}\to \frac{(\gamma+1)^2}{4k_2\gamma}=\frac{k_2}{k_2^2-k_1^2}$$
as $r\to 1$. Then, it is clear that
$$\lim_{t\to 0}\frac{x(t)-\frac\lambda 2\pi i}{t} = k_1^2+k_2^2.$$
\end{proof}

\section{Vertical movement of poles}\label{vert}

In this section we study the vertical motion of the singularities of the two--soliton solutions $u^\pm$ of (\ref{mkdv}).  As seen in Section \ref{formula}, this comes down to studying the zeroes of $F^\pm$ given by
\re{defF+}--\re{defF-}.  Recall the solution $u^+$ represents two interacting solitons of the same sign, while $u^-$ represents two interacting solitons of opposite sign.  Observe that 
\begin{align}
F^+(x,t) = F^+_1(x,t) F^+_2(x,t)\\
F^-(x,t)  = F^-_1(x,t) F^-_2(x,t)
\end{align}
where
\begin{align}
F^+_1
= 1 + i\gamma f_1 + i\gamma f_2 - f_1f_2\\
F^+_2 = 1 - i\gamma f_1 - i\gamma f_2 - f_1f_2\\
F^-_1
= 1 + i\gamma f_1 - i\gamma f_2 + f_1f_2\\
F^-_2 =1 - i\gamma f_1 + i\gamma f_2 + f_1f_2
\end{align}
and $f_1$ and $f_2$ are as in (\ref{deffj}).
Note the similarity in form between the above functions
and the function $F$ defined by formula (2.13) in \cite{BW2}.

Since $F^+_1(x,t) = 0$ if and only if $F^+_2(\overline x,t) = 0$, to study the
zeros of $F^+$, it suffices to study the zeros of $F^+_1$.  Similarly, since
$F^-_1(x,t) = 0$ if and only if $F^-_2(\overline x,t) = 0$, to study the
zeros of $F^-$, it suffices to study the zeros of $F^-_1$.

As before, it is sometimes necessary to distinguish between the cases where
 $k_1$ and $k_2$ are commensurable and the cases where they
 are not. In the former case,  denote by
 $p_1$ and $p_2$ the relatively prime positive integers such that
(\ref{defCom}) holds, and let $\lambda$ be defined as in \re{defLam}.
Hence, if $k_1$ and $k_2$ are commensurable, the functions $f_1$ and $f_2$ 
are $2\pi \lambda i$ periodic in $x$.  
Thus, in the commensurable case, all four of the functions $F_1^\pm, F_2^\pm$,
must also be  $2\pi \lambda i$ periodic in $x$.

\begin{proposition} The zeroes of $F^+$ and $F^-$ in the complex plane lie 
off the real axis.    Moreover, in the commensurable case, if either $F^+(x,t) = 0$ or $F^-(x,t) = 0$,
then $\Im x \not = 2m\pi \lambda$ for all $m \in \Z$.
\end{proposition}

\begin{proof}  Suppose $x\in\R$ and
$F^+(x,t) = 0$.  It follows that $F^+_1(x,t) = F^+_2(x,t) = 0$.  Taking real and imaginary parts 
produces the coupled system
\begin{align*}
1  &= f_1f_2,\\
f_1 &= -f_2,
\end{align*}
from which it follows that $f_1(x,t)^2 = -1$, which is impossible.
If $x\in\R$ and
$F^-(x,t) = 0$,  it similarly follows that
\begin{align*}
1  &= -f_1f_2,\\
f_1 &= f_2,
\end{align*}
from which one observes that $f_1(x,t)^2 = -1$, which is again impossible.

The last statement follows from the $2\pi \lambda i$ periodicity which holds 
in the commensurable case.
\end{proof}

If $k_1$ and $k_2$ are commensurable, further information
about the location of the zeros of $F^\pm$ can be obtained.

\begin{proposition} Suppose that $k_1$ and $k_2$ are commensurable. 
If either $F^+(x,t) = 0$ or $F^-(x,t) = 0$,
then $\Im x \not = (2m+1)\pi \lambda$ for any $m \in \Z$.
\end{proposition}

\begin{proof}
Because of $2\pi \lambda i$--periodicity, it suffices to consider $m = 0$.
Suppose $\Im x = \pi\lambda$ and $F^+_1(x,t) = 0$.
Since $x - \overline x = 2\pi\lambda i$, it follows from 
$2\pi \lambda i$-periodicity that $F^+_1( \overline x,t) = 0$.  But this says precisely that
$F^+_2(x,t) = 0$.  Adding and subtracting the two equations,
$F^+_1(x,t) = 0$ and $F^+_2(x,t) = 0$ gives
\begin{align*}
1  &= f_1f_2,\\
f_1 &= -f_2,
\end{align*}
from which it is concluded that $f_1(x,t)^2 = -1$, i.e. $f_1(x,t) = \pm i$ and $f_2(x,t) = \mp i$.  
The latter system together with the definition (\ref{deffj}) of the $f_j$'s implies
\begin{align*}
\Re x = k_1^2 t,\\
\Re x = k_2^2 t,
\end{align*}
and so $t = \Re x = 0$.  Thus $x = \pi\lambda i$.  However,
$$
f_1( \pi\lambda i,0) = \exp(-k_1\pi\lambda i) = \exp(-p_1\pi i) = \pm 1,
$$
since $p_1$ is an integer.  This contradiction shows that
$F^+_1(x,t)$ can not be zero if $\Im x = \pi\lambda$.  A similar
argument applies to the other functions $F^+_2$, $F^-_1$ and
$F^-_2$.
\end{proof}

To show that most of the poles of the two--soliton solution of \re{mkdv} move vertically, 
we need to reproduce calculations which are similar to
those found in Section 3 of \cite{BW2}.  Unfortunately, it seems that the calculations in 
\cite{BW2} do not directly apply to the present situation in all cases.  In an effort at 
economy, we only prove here that, with certain very specific exceptions, the poles in 
the two--soliton solution always feature vertical movement.   
Consequently, we do not need to reproduce the entirety of
Section 3 of \cite{BW2}.  Nonetheless, we cannot avoid certain, somewhat tedious 
calculations,  closely modeled on those in \cite{BW2}.

The following notation
$$\alpha = -\Im x,$$
$$A_1 = e^{-k_1\Re x + k_1^3 t},$$
$$A_2 = e^{-k_2\Re x + k_2^3 t}, $$
is taken from \cite{BW2}.  
The function $F_1^+$ may be rewritten in this notation, {\it viz.}
\begin{equation}\label{F+1alpha}
F^+_1(x,t) = 1 + i\gamma A_1 e^{ik_1\alpha} + i\gamma A_2 e^{ik_2\alpha} - A_1A_2e^{i(k_1+k_2)\alpha}.
\end{equation}

Attention is first focussed on the solution $u^+$.  To investigate possible vertical 
movement of poles of $u^+$, we examine more closely the zeros of $F^+_1(x,t)$.

\begin{proposition}\label{bddimpt}
Suppose $F^+_1(x,t) = 0$.  Then, 
$\cos k_1 \alpha = 0$ if and only if $\cos k_2 \alpha = 0$.
Moreover, the relation
\begin{equation}
\label{cosk12}\Big(A_2 - \frac{1}{A_2}\Big)\cos k_1\alpha
+ \Big(A_1 - \frac{1}{A_1}\Big)\cos k_2\alpha = 0
\end{equation}
always holds.
In case $\cos k_1 \alpha \not = 0$ and $\cos k_2 \alpha \not = 0$,
then  $A_1 = 1$ if and only if $A_2 = 1$, and this can
only happen if $t = 0$ and $\Re x = 0$.

\end{proposition}

\begin{proof}
First, multiply the equation $F^+_1(x,t) = 0$
 by $1 - i\gamma A_1 e^{-ik_1\alpha}$, express this using the representation (\ref{F+1alpha}) 
and take the imaginary part of the result. This leads to the formula
\begin{equation}
\label{cosk2}\Big(A_1 + \frac{1}{A_1}\Big)\cos k_2\alpha  = \frac{1}{\gamma}\sin (k_2 + k_1)\alpha
- \gamma \sin (k_2 - k_1)\alpha.
\end{equation}
In the same way, multiplying by $1 - i\gamma A_2 e^{-ik_2\alpha}$ and subsequently
 extracting the imaginary part yields 
\begin{equation}
\label{cosk1}\Big(A_2 + \frac{1}{A_2}\Big)\cos k_1\alpha  = \frac{1}{\gamma}\sin (k_2 + k_1)\alpha
+ \gamma \sin (k_2 - k_1)\alpha.
\end{equation}
If $\cos k_2\alpha = 0$, so that in particular $\sin k_2\alpha  = \pm 1$, it follows from (\ref{cosk2}),
using the formulas for the sine of the sum and difference, that $\cos k_1\alpha = 0$.
In the same way, using (\ref{cosk1}), if $\cos k_1\alpha = 0$, then $\cos k_2\alpha = 0$.  This proves the first assertion in the proposition.

For future reference, notice that it follows by subtracting (\ref{cosk2}) from  (\ref{cosk1})
that
\begin{equation}
\label{sink21}\gamma \sin (k_2 - k_1)\alpha =  \frac{1}{2} \Big(A_2 + \frac{1}{A_2}\Big)\cos k_1\alpha - \frac{1}{2}\Big(A_1 + \frac{1}{A_1}\Big)\cos k_2\alpha.
\end{equation}

Taking the
real and imaginary parts of (\ref{F+1alpha}) and setting them equal to zero yields
\begin{align}
1 - \gamma A_1 \sin k_1\alpha - \gamma A_2\sin k_2\alpha
- A_1A_2 \cos(k_1+k_2)\alpha &= 0,  \\
\label{impart} \gamma A_1 \cos k_1\alpha + \gamma A_2\cos k_2\alpha
- A_1A_2 \sin(k_1+k_2)\alpha &= 0.
\end{align}
Multiply the first equation above by $\sin(k_1+k_2)\alpha$, the second
equation by $\cos(k_1+k_2)\alpha$ and form the difference of the outcomes.  
The formula
\begin{equation}
\sin(k_1+k_2)\alpha = \gamma A_2 \cos k_1 \alpha +  \gamma A_1 \cos k_2 \alpha
\end{equation}
emerges from these machinations.  But, (\ref{impart}) implies
\begin{equation}
\sin(k_1+k_2)\alpha = \frac{\gamma}{A_2} \cos k_1 \alpha +  \frac{\gamma}{ A_1} \cos k_2 \alpha.
\end{equation}
The last two equations taken together imply (\ref{cosk12})

This proves the second assertion of the proposition.
Finally, it is clear that $A_1 = A_2 = 1$ if and only if
$t = 0$ and $\Re x = 0$.
\end{proof}

\begin{proposition}\label{alph}
Suppose that $\cos k_1 \alpha = \cos k_2 \alpha = 0$.  It follows that
$k_1$ and $k_2$ are commensurable, that
$p_1$ and $p_2$ are both odd and that $\alpha$ is an odd multiple of
$\pi\lambda/2$.
\end{proposition}

\begin{proof}
If $\cos k_1 \alpha = \cos k_2 \alpha = 0$, then there exist integers $m$ and $n$ such that
\begin{align*}
k_1 \alpha  = (2m + 1)\frac{\pi}{2},  \\
k_2 \alpha  = (2n + 1)\frac{\pi}{2}.
\end{align*}
Thus it is clear that $k_1$ and $k_2$ are commensurable and
\begin{equation*}
\frac{p_1}{p_2} = \frac{k_1}{k_2} = \frac{2m + 1}{2n + 1},
\end{equation*}
whence
\begin{equation}
\label{podd}
p_1(2n + 1) = p_2(2m + 1).
\end{equation}
It follows that $p_2 - p_1$ is even, and since they are also relatively
prime, they both must be odd.
Next, equation (\ref{podd}) implies that $p_1 = (2m+1)/c$, where 
$$c = \gcd(2m+1,2n+1).
$$
Consequently, 
\begin{equation*}
\alpha = \frac{2m + 1}{k_1}\frac{\pi}{2} = \frac{2m + 1}{p_1}\frac{\pi\lambda}{2}  = c\frac{\pi\lambda}{2},
\end{equation*}
which concludes the proof since $c$ is necessarily odd.
\end{proof}

\begin{proposition}
Let $x(t)$ be a smooth curve such that $F^+_1(x(t),t) = 0$ and
$\partial_x F^+_1(x(t),t) \not = 0$.  It follows that $\Im x'(t)$ has the same
sign as
\begin{equation*}
\Big(A_1 - \frac{1}{A_1}\Big)\cos k_2\alpha
= - \Big(A_2 - \frac{1}{A_2}\Big)\cos k_1\alpha.
\end{equation*}
\end{proposition}

\begin{proof} Let $x(t)$ be a smooth curve such that $F^+_1(x(t),t) = 0$ and
$\partial_xF^+_1(x(t),t) \not = 0$.
Implicit differentiation provides the relation
\begin{align*}
x'(t) &= -\frac{\partial_tF^+_1(x(t),t)}{ \partial_x F^+_1(x(t),t)}    \\[.3cm]  
&= \frac{k_1^3i\gamma A_1 e^{ik_1\alpha} + k_2^3i \gamma A_2 e^{ik_2\alpha}
- (k_1^3 + k_2^3) A_1A_2 e^{i(k_1 + k_2)\alpha}}
{ k_1 i\gamma A_1 e^{ik_1\alpha} + k_2 i\gamma A_2 e^{ik_2\alpha}
- (k_1 + k_2) A_1A_2 e^{i(k_1 + k_2)\alpha}}\\[.4cm]
&= \frac{k_1^3\gamma A_1 e^{-ik_2\alpha} + k_2^3 \gamma A_2 e^{-ik_1\alpha}
+ i(k_1^3 + k_2^3) A_1A_2}
 {k_1 \gamma A_1 e^{-ik_2\alpha} + k_2 \gamma A_2 e^{-ik_1\alpha}
+i (k_1 + k_2) A_1A_2}
\end{align*}
\addtolength{\abovedisplayskip}{-12pt}
\begin{align*} \! \! = 
  \frac{1}{\big| k_1 \gamma A_1 e^{-ik_2\alpha} + k_2 \gamma A_2 e^{-ik_1\alpha}
+ i(k_1 + k_2) A_1A_2\big|^2}  \Big[
 \big(k_1^3 \gamma A_1 e^{-ik_2\alpha} + k_2^3 \gamma A_2 e^{-ik_1\alpha}  \\ + i(k_1^3 + k_2^3) A_1A_2 \big)
\big(k_1 \gamma A_1 e^{ik_2\alpha} + k_2 \gamma A_2 e^{ik_1\alpha}
- i(k_1 + k_2) A_1A_2\big)  \Big].
\end{align*}
The imaginary part of the numerator in this last expression is equal to
\begin{equation*}
\addtolength{\abovedisplayskip}{12pt}
k_1 k_2(k_2^2 - k_1^2)\gamma^2 A_1A_2 \sin (k_2 - k_1)\alpha
 + k_1k_2(k_1 + k_2)^2A_1A_2(A_1\cos k_2\alpha -  A_2\cos k_1\alpha),
\end{equation*}
which is a positive mulitple of, and so has the same sign as,
\addtolength{\abovedisplayskip}{12pt}
\begin{equation}\label{signxprime}
\gamma  \sin (k_2 - k_1)\alpha
 + A_1\cos k_2\alpha -  A_2\cos k_1\alpha.
\end{equation}
Finally, formulas (\ref{sink21}) and  (\ref{cosk12}) reveal that the expression in (\ref{signxprime}) is equal to
\begin{equation*}
A_1 \cos k_2\alpha -  \frac{1}{A_1}\cos k_2\alpha.
\end{equation*}
\end{proof}

\begin{corollary}Suppose either that $k_1$ and $k_2$ are not commensurable, or that they are commensurable with $p_1$ and $p_2$ having opposite parity.
Let $x(t)$ be a smooth curve such that $F^+_1(x(t),t) = 0$. If either $t \not = 0$ or if $\Re x \not = 0$, then $\Im x'(t) \not = 0$.
 \end{corollary}

 \begin{proof}
We know from Lemma \ref{lem-multzeros} that
the roots of $F^+(\cdot,t)$ are all simple unless
$k_1$ and $k_2$ are commensurable and $p_1$ and $p_2$ are both odd.
It follows that the same is true for $F^+_1(\cdot,t)$.
Thus, it must be the case that $\partial_xF^+_1(x(t),t) \not = 0$.
The last three propositions can now be brought to bear to establish the claim.
 \end{proof}

To recapitulate, it has been shown, except in the commensurable case with
$p_1$ and $p_2$ both odd, that if either
$t \not = 0$ or if $\Re x \not = 0$, then $\Im x'(t) \not = 0$,
so the poles of the solution $u^+$ are {\it always} moving
vertically.  The same is true in the case that $p_1$ and $p_2$ are both odd,
so long as the imaginary part of the pole is not an odd
multiple of $\pi\lambda/2$.
 The same conclusions are true about the poles of $u^-$.
The calculations leading to this conclusion, while not the
same as those for $u^+$, are analogous enough that we 
pass over the details.   Interestingly, 
in the commensurable case when $p_1$ and $p_2$ have
opposite parity, it turns out that the poles of $u^-$
and those of $u^+$ bear a very simple relationship to 
one another.

\begin{proposition}  Suppose that $p_1$ and $p_2$ have
opposite parity (one odd, one even).  It follows that there
exists $\theta \in \R$ such that
\begin{equation}\label{translate}
F^+(x - i\theta,t) = F^-(x,t)
\end{equation}
for all $x \in \C$ and $t \in \R$.  In other words, the
poles of $u^-$ are precisely given by a vertical translation
of the poles of $u^+$.
\end{proposition}

\begin{proof} Suppose $\theta \in \R$ is such that
\begin{align}
\label{f1ef1}f_1(x - i\theta,t) &= f_1(x,t),\\
\label{f2emf2}f_2(x - i\theta,t) &= -f_2(x,t),
\end{align}
for all $x \in \C$ and $t \in \R$.  It would then follow that
(\ref{translate}) holds
for all $x \in \C$ and $t \in \R$.  The same would be true if we had
instead
\begin{align}
\label{f1emf1}f_1(x - i\theta,t) &= -f_1(x,t),\\
\label{f2ef2}f_2(x - i\theta,t) &= f_2(x,t).
\end{align}
For (\ref{f1ef1}) and  (\ref{f2emf2}) to be valid, it is necessary and sufficient that
\begin{align}
\exp(ik_1\theta) &= \exp(ip_1\theta/\lambda) = 1,\\
\exp(ik_2\theta) &= \exp(ip_2\theta/\lambda) = -1.
\end{align}
For these latter conditions to hold, it is necessary for there to be 
two integers $m$ and $n$ such that
\begin{align}
p_1\theta/\lambda &= 2m\pi,\\
p_2\theta/\lambda &= (2n+1)\pi,
\end{align}
or, what is the same, 
\begin{equation}
\frac{\theta}{\lambda\pi} = \frac{2m}{p_1} = \frac{2n+1}{p_2}.
\end{equation}
If $p_1$ is even and $p_2$ is odd, it is clear that one may choose appropriate values of $\theta, m$ 
and $n$ so that the last equation holds

In the opposite case, if $p_1$ is odd and $p_2$ is even, a similar argument
shows that there exists  $\theta, m$and $n$ such that (\ref{f1emf1}) and
(\ref{f2ef2}) hold.
\end{proof}

It remains to consider the case where $k_1$ and $k_2$
are commensurable, with $p_1$  and $p_2$ odd.
To establish the existence of poles with non--trivial vertical
movement, it suffices by Proposition~\ref{alph} to establish the
existence of poles whose imaginary parts are not an odd multiple
of $\pi\lambda/2$ (with either $t \not = 0$ or with non--zero real part).

As shown in \cite{BW2}, for the two--soliton solution of the KdV--equation
in the commensurable case, there are always poles 
 with vertical movement, and at least one pole moving
precisely horizontally. The proof of this fact requires the full force
of the delicate and technical analysis in \cite{BW2}.  We would like to
avoid such technical calculations in this paper to the extent possible.

Thus, for the modified KdV--equation (\ref{mkdv}),  we should not expect that all the poles will be moving vertically in the remaining cases.  The fact that in the case of opposite parity,
all the poles move vertically (except if $t = 0$ and $\Re x = 0$) is already
an interesting difference in behavior between the two equations.

It turns out that in the commensurable case, with $p_1$  and $p_2$ both odd,
the movement of the poles of $u^+$ can in fact be reduced to the movement
of the poles of the 2--soliton solution of the KdV--equation if $p_2 - p_1 \in 4\N$.
The same is true for $u^-$ if $p_2 + p_1 \in 4\N$.
Here are the precise statements.

\begin{proposition}  Suppose that $p_1$ and $p_2$ are both odd.
If $p_2 - p_1 \in 4\N$ then there
exist $\theta_1,\theta_2 \in \R$ such that
\begin{equation}\label{translate1}
F^+_1(x - i\theta_1,t) =  1 + \gamma f_1(x,t) + \gamma f_2(x,t) + f_1(x,t)f_2(x,t)
\end{equation}
and
\begin{equation}\label{translate2}
F^+_2(x - i\theta_2,t) =  1 + \gamma f_1(x,t) + \gamma f_2(x,t) + f_1(x,t)f_2(x,t)
\end{equation}
for all $x \in \C$ and $t \in \R$.
\end{proposition}

\begin{proof} Suppose $\theta_1 \in \R$ is such that
\begin{align}
\label{f1emif1}f_1(x - i\theta_1,t) &=- if_1(x,t),\\
\label{f2emif2}f_2(x - i\theta_1,t) &= -if_2(x,t),
\end{align}
for all $x \in \C$ and $t \in \R$.  It would follow that
(\ref{translate1}) holds 
for all $x \in \C$ and $t \in \R$.
In addition, if we have
\begin{align}
\label{f1eif1}f_1(x - i\theta_2,t) &=if_1(x,t),\\
\label{f2eif2}f_2(x - i\theta_2,t) &= if_2(x,t),
\end{align}
then (\ref{translate2}) would be true.

For the system (\ref{f1emif1})--(\ref{f2emif2}) to be valid, it is necessary and sufficient that
\begin{align}
\exp(ik_1\theta) &= \exp(ip_1\theta/\lambda) = -i,\\
\exp(ik_2\theta) &= \exp(ip_2\theta/\lambda) = -i.
\end{align}
For this, we need to find two integers $m$ and $n$ such that
\begin{align}
\frac{p_1\theta}{\lambda} &= (4m-1)\frac{\pi}{2},\\
\frac{p_2\theta}{\lambda} &= (4n-1)\frac{\pi}{2},
\end{align}
which is the same as asking for two integers $m$ and $n$ such that 
\begin{equation}
\frac{2\theta}{\lambda\pi} = \frac{4m-1}{p_1} = \frac{4n-1}{p_2}.
\end{equation}
For the latter to hold true, it must be the case that
\begin{equation}
\frac{4(p_2m - p_1n)}{p_2 - p_1} = 1.
\end{equation}
Since $p_1$ and $p_2$ are relatively prime, there exist integers
$r$ and $s$ such that
$$
rp_2 + sp_1= 1;
$$
hence, simply take  
\begin{align*}
m &= \frac{(p_2-p_1)r}{4},\\
n &= -\frac{(p_2-p_1)s}{4}.
\end{align*}

The proof for $F^+_2$ is similar, but with $4m+1$ replacing
$4m-1$ and $4n+1$ replacing
$4n-1$.
\end{proof}

\begin{proposition}  Suppose that $p_1$ and $p_2$ are both odd.
If $p_2 + p_1 \in 4\N$ then there
exist $\theta_1,\theta_2 \in \R$ such that
\begin{equation}\label{translate3}
F^-_1(x - i\theta_1,t) =  1 + \gamma f_1(x,t) + \gamma f_2(x,t) + f_1(x,t)f_2(x,t)
\end{equation}
and
\begin{equation}\label{translate4}
F^-_2(x - i\theta_2,t) =  1 + \gamma f_1(x,t) + \gamma f_2(x,t) + f_1(x,t)f_2(x,t)
\end{equation}
for all $x \in \C$ and $t \in \R$.
\end{proposition}

The expression on the right side of the four formulas
(\ref{translate1}), (\ref{translate2}), (\ref{translate3}), and (\ref{translate4}),
i.e.
\begin{equation}
1 + \gamma f_1(x,t) + \gamma f_2(x,t) + f_1(x,t)f_2(x,t),
\end{equation}
is exactly the function $F$ in formula (2.13) of \cite{BW2}
whose zeros correspond to the poles of the 2--soliton solution of the 
KdV--equation.
Thus, we may use the results of \cite{BW2} to describe the behavior
of the poles of $u^\pm$ in the  cases under consideration.

More precisely, we may now affirm that in the case where
  $p_1$ and $p_2$ are both odd and
 $p_2 - p_1 \in 4\N$, the solution $u^+$ has two poles moving
horizontally on the line $\Im x = \pi\lambda/2$ and also on the line
$\Im x = -\pi\lambda/2$.  There are $2(p_1 + p_2 - 2)$ other poles
with imaginary part between $\pm\pi\lambda$, and they will move vertically
for all $t \not = 0$.  The same will
be true for $u^-$ in the case $p_1$ and $p_2$ are both odd if $p_2 + p_1 \in 4\N$.  These configurations will repeat
with $2\pi\lambda i$ periodicity, so that horizontally moving poles
are found with imaginary part equal to every odd multiple of
$\pi\lambda/2$.

The last case is  one where four poles meet
at $t = 0$,  as described in the previous section, at  poles 
on the imaginary axis whose imaginary part
is an odd multiple of $\pi\lambda/2$.  At each of these points,
the analysis shows that (at least) two of these poles have vertical
movement.  We refrain from going into the detailed considerations 
needed to establish that all the poles move vertically for all $t \not = 0$,
excepting the two poles approaching the singular points horizontally.  
In any event, we have shown in this case the existence
of poles with vertical movement.

\section{Finite time blowup of solutions}

The results of the previous section immediately imply that
there exist complex--valued solutions to (\ref{mkdv}) on $\R$
which blow up in finite time.  Indeed, fix $\alpha \in \R$,
and let $u$ be given by
\begin{equation}\label{comsol}
u(x,t)  = u^\pm(x - i\alpha,t),
\end{equation}
where we may use either of the two soliton solutions
$u^+$ or $u^-$.  As long as the set $\{(x - i\alpha,t_0): x \in \R\}$ does not 
contain a pole
of $u^\pm$, $u$ is a smooth, complex--valued solution of (\ref{mkdv}) for
$t$ in a neighborhood of $t_0$.     Since at any given time, 
 the imaginary parts of the collection 
of all poles of $u^+$ or $u^-$  form
a discrete set, and since both $u^+$ and $u^-$ have poles whose
imaginary parts move continuously in time, there exist $\alpha$ and $t_0$
such that $(x - i\alpha,t_0)$ is a pole of $u^+$,  say, but $(x - i\alpha,t)$ is not a pole of $u^+$ if
$t$ is sufficiently close to, but not equal to, $t_0$.
It follows that the solution $u$ defined in  (\ref{comsol}) with this choice of $\alpha$ is a
regular complex--valued solution of (\ref{mkdv}) on $\R$, decaying
exponentially to zero as $x \to \pm \infty$ for $t$ close
to $t_0$, but which is singular at $t = t_0$.  In other words, the solution
blows up in finite time.

It is interesting to note that this result of singularity formation for complex--valued
solutions of the mKdV--equation can be interpreted
as a  blow--up result for real--valued solutions of a system of dispersive
equations.  Let $u$ be a solution of (\ref{mkdv}) and let
$r = \Re u$ and $s = \Im u$.  It follows that $r$ and $s$ satisfy the
real--valued system
\begin{align}
r_t+r_{xxx}+6(r^2 - s^2)r_x - 2rss_x = 0, \\
s_t+s_{xxx} + 2rsr_x+6(r^2 - s^2)s_x = 0.
\end{align}
Thus, we have shown that this coupled, dispersive system admits
real--valued solutions (exponentially decaying in space) which
blow up in finite time.


\section{Some formal calculations}\label{maple}
(All the computations in this section have been carried out using MAPLE.)
As noted at the end of Section \ref{formula}, the interaction
time for the two--soliton solutions $u^\pm$, given in
\re{uplus} and \re{uminus}, is $t = 0$, the center
of the interaction is $x = 0$ and at $t=0$, the solution 
is even in $x$.   The explicit formulas
for $u^\pm$ allow us to observe and calculate certain
aspects of these solutions at the moment of interaction.
In particular, it is interesting to know whether
there is a single maxima during the interaction or not, and it 
is likewise interesting to know the speed of the
two solitons at the moment
of interaction.

In the case of $u^+$, one observes that
that the solution has one centered maximum
if the ratio $k_2/k_1$ is large enough (bigger than
around 2.6), and two symmetrically located maxima
for smaller values of $k_2/k_1$.
Partial confirmation of this can be obtained 
by computing $u^+_{xx}(0,0)$.  Since, by symmetry,
$u^+_x(0,0) = 0$, the sign of the second derivative
will tell us if it is a local maximum or a local minimum.
A MAPLE--implemented calculation shows that 
\begin{equation*}
u^+_{xx}(0,0) = -(k_2 - k_1)(k_1^2 - 3k_1k_2 + k_2^2).
\end{equation*}
Therefore, if $1 < k_2/k_1 < (3 + \sqrt5)/2$, then $u^+_{xx}(0,0) > 0$,
which means that $u^+(\cdot,0)$ has a local minimum at $x = 0$.
Thus, there are (at least) two maxima at the moment of interaction.
If $k_2/k_1 > (3 + \sqrt5)/2$, then $u^+(\cdot,0)$ has a local maximum at $x = 0$, which is consistent with there being a single maximum
at the moment of interaction.

Similarly, we calculate that
\begin{equation*}
u^-_{xx}(0,0) = -(k_2 + k_1)(k_1^2 + 3k_1k_2 + k_2^2),
\end{equation*}
which means $u^-(\cdot,0)$ has a local maximum at $x = 0$
for all values $0 < k_1 < k_2$, which is consistent with
the graphical observations of the solution itself.

If we are interested in the ``speed" of the two--soliton
solution at the moment of interaction, one idea
is to calculate the speed of
the maximum  or the local minimum  of the solution as $(x,t)$ approaches $(0,0)$. 
 One might
argue that the movement of the  maximum is some kind of speed.  It must be 
acknowledged that the interpretation is less clear when one is tracking a 
local minimum.
To calculate this speed, suppose $y(t) = y_\pm(t)$ is a real--valued
curve such that $u^\pm_x(y(t),t) = 0$, i.e. $y_\pm(t)$
is always at an extemal point of the solution $u = u^\pm$.
Suppose also  that $y_\pm(0) = 0$, which is to say the curve is at the 
interaction center at the
interaction time $t = 0$.
Differentiating with respect to $t$, we see that
 $u_{xx}(y(t),t)y'(t) + u_{xt}(y(t),t) = 0$, or
 \begin{equation*}
 y'(t) = -\frac{u_{xt}(y(t),t)}{u_{xx}(y(t),t)}.
 \end{equation*}
 This gives, in turn,
  \begin{equation}
 y'(0) = -\frac{u_{xt}(0,0)}{u_{xx}(0,0)},
 \end{equation}
 which might be thought of as representing the speed of the two--soliton
 solution at the moment of interaction of the two solitons.
 The value of $y'(0)$ can be calculated explicitly
 from \re{uplus} and \re{uminus}.
 The results are as follows, for both $u^+$ and $u^-$,
 \begin{align*}
 y_+'(0) = \frac{k_1^4 - 3k_1^3k_2 + 3k_1^2k_2^2 - 3k_1k_2^3 + k_2^4}
 {k_1^2 - 3k_1k_2 + k_2^2},\\
  y_-'(0) = \frac{k_1^4 + 3k_1^3k_2 + 3k_1^2k_2^2 + 3k_1k_2^3 + k_2^4}
 {k_1^2 + 3k_1k_2 + k_2^2}.
 \end{align*}
For $u^-$, where there is always a maximum at $x = 0$ at the moment  of interaction, the  maximum is moving with a
 posistive speed.  What that speed represents is not entirely clear.
 In the case of $u^+$, where the midpoint is a local maximum
 only if $k_2/k_1 > (3 + \sqrt5)/2$, we see that for these values
 of $0 < k_1 < k_2$, we have indeed $  y_+'(0) > 0$.  On the other hand,
 $y_+'(0) < 0$ for at least some values of $0 < k_1 < k_2$ with
 $k_2/k_1 < (3 + \sqrt5)/2$.  (The lower bound on  $k_2/k_1$ for which
 this speed is negative is around 2.15.)  This negative speed
 represents the speed of the local minimum, between the two maxima.
 It is curious that in some cases this minimum
 is moving backwards.

 \bigskip
 It is also interesting to do this for the two--soliton
 solution of the KdV--equation, a calculation which was not carried out in
 \cite{BW2}.  In this case, we obtain
 \begin{equation*}
 y'(0) = \frac{k_1^4  + 2k_1^2k_2^2  - k_2^4}
 {3k_1^2 - k_2^2}.
 \end{equation*}
 As is known, the two--soliton solution of the KdV--equation has one
 maximum at the interaction time if $k_2/k_1 > \sqrt3$
 and two maxima, symmetrically located about the interaction center,  if $1 < k_2/k_1 < \sqrt3$.
 Here we see that $y'(0) < 0$ for $\sqrt{1 + \sqrt2} < k_2/k_1 < \sqrt3$.
 In these case, of course, it is the minimum which is moving backward.

\end{document}